\documentclass[psamsfonts]{amsart}
\input{diagrams.tex}
\diagramstyle[scriptlabels,height=8mm,width=8mm]
\newcommand{\query}[1]%
{\mbox{}\marginpar{\raggedright\hspace{0pt}{\small\em #1}}}%
\theoremstyle{plain}

\newtheorem{thm}{Theorem}[section]
\newtheorem{cor}[thm]{Corollary}
\newtheorem{lem}[thm]{Lemma}
\newtheorem{prop}[thm]{Proposition}

\theoremstyle{definition}
\newtheorem{defi}[thm]{Definition}
\newtheorem{conj}[thm]{Conjecture}
\newtheorem{conv}[thm]{Convention}
\newtheorem{nota}[thm]{Notation}
\newtheorem{rem}[thm]{Remark}
\newtheorem{exa}[thm]{Example}
\newtheorem{sit}[thm]{}

\newcommand{\brem}{\begin{rem}}
\newcommand{\erem}{\end{rem}}
\newcommand{\bexa}{\begin{exa}}
\newcommand{\eexa}{\end{exa}}
\newcommand{\bdefi}{\begin{defi}}
\newcommand{\edefi}{\end{defi}}
\newcommand{\bcor}{\begin{cor}}
\newcommand{\ecor}{\end{cor}}
\newcommand{\blem}{\begin{lem}}
\newcommand{\elem}{\end{lem}}
\newcommand{\bconv}{\begin{conv}}
\newcommand{\econv}{\end{conv}}
\newcommand{\bconj}{\begin{conj}}
\newcommand{\econj}{\end{conj}}
\newcommand{\bprop}{\begin{prop}}
\newcommand{\eprop}{\end{prop}}
\newcommand{\bthm}{\begin{thm}}
\newcommand{\ethm}{\end{thm}}
\newcommand{\bnota}{\begin{nota}}
\newcommand{\enota}{\end{nota}}
\newcommand{\bsit}{\begin{sit}}
\newcommand{\esit}{\end{sit}}

\def\fp{{\mathfrak p}}

\def\cH{{\mathcal H}}
\def\cL{{\mathcal L}}
\def\cO{{\mathcal O}}

\def\gL{{\Lambda}}
\def\gD{{\Delta}}

\def\tD{{\widetilde D}}
\def\hA{{\widehat A}}

\newcommand{\A}{{\mathbb A}}
\newcommand{\G}{{\mathbb G}}

\newcommand{\C}{{\mathbb C}}
\newcommand{\Q}{{\mathbb Q}}
\newcommand{\Z}{{\mathbb Z}}
\newcommand{\N}{{\mathbb N}}

\def\lto{\longrightarrow}
\def\hto{\hookrightarrow}
\def\ul{\underline}

\newcommand{\Ext}{\operatorname{Ext}}
\newcommand{\Sing}{\operatorname{Sing}}
\newcommand{\Spec}{\operatorname{Spec}}

\newcommand{\red}{{\operatorname{red}}}
\newcommand{\Proj}{\operatorname{Proj}}

\newcommand{\pr}{{\operatorname{pr}}}
\newcommand{\codim}{{\operatorname{codim}}}

\newcommand{\la}{\label}

\newcommand{\be}{\begin{eqnarray}}
\newcommand{\ee}{\end{eqnarray}}
\newcommand{\no}{\noindent}

\title{Log-canonical forms  and log canonical
singularities}

\author{Hubert Flenner}
\address{Fakult\"at f\"ur Mathematik,
Ruhr Universit\"at Bochum,
Geb.\ NA 2/72,
Universit\"ats\-str.\ 150,
44780 Bochum, Germany}
\email{Hubert.Flenner@ruhr-uni-bochum.de}

\author{Mikhail Zaidenberg}
\address{Universit\'e
Grenoble I, Institut Fourier,
UMR 5582 CNRS-UJF, BP 74,
38402 St.\ Martin
d'H\`eres c\'edex, France}
\email{zaidenbe@ujf-grenoble.fr}

\thanks{This work was supported by stays of the
first  author at the University of Grenoble
and by both authors at the Max
Planck Institute of Mathematics at Bonn;  we would
like to thank these institutions for their
hospitality.\\
\mbox{\hspace{11pt}}{\it 1991 Mathematics
Subject Classification}:
14B05, 14B15, 14E15, 14J17, 32S20, 32S25.\\
\mbox{\hspace{11pt}}{\it Key words}:
rational
singularity, log canonical singularity,
log terminal singularity,
quotient singularity, plurigenus,
quasihomogeneous singularity, $\C^*$-action,
graded algebra, Kodaira dimension}

\date{\today}

\begin{document}
\maketitle
\centerline{\it Dedicated to
H.\ Grauert on occasion of his
70's birthday}

\begin{abstract}
For a normal subvariety $V$ of $\C^n$ with
a good $\C^*$-action we give a
simple characterization for when
it has only log canonical, log terminal or rational
singularities.  Moreover  we are able to
give formulas for the plurigenera
of isolated singular points of such
varieties and of the logarithmic
Kodaira dimension of $V\backslash
\{0\}$.  For this purpose
we introduce sheaves
of $m$-canonical and $L^{2,m}$-canonical
forms on normal complex
spaces. For the case of affine
varieties with good
$\C^*$-action we give an
explicit formula for these sheaves in terms
of the grading of the dualizing
sheaf and its tensor powers.
\end{abstract}

\tableofcontents

\section*{Introduction}

Let $X$ be a normal complex space and $D$
a reduced Weil divisor on $X$.
In this paper we will associate
to the pair $(X,D)$ two sheaves
$\cL_{X,D}^m$ and $\cL_{X,D}^{2,m}$,
which we call the sheaves of
logarithmic $m$-, resp.\
$L^{2}$-$m$-canonical forms.  The
construction is in brief
as follows: let $\pi:X'\to X$ be a
resolution of singularities
such that $D':=\pi^{-1}(D\cup \Sing
X)_\red$ is an SNC
(=simple normal crossing) divisor.
The sheaves
$$
\cL_{X,D}^m:=\pi_*(\cO_{X'}(m(K_{X'}+D')))
\quad
\mbox{and}\quad
\cL_{X,D}^{2,m}:=
\pi_*(\cO_{X'}(mK_{X'}+(m-1)D'))
$$
are then independent of
the choice of resolution.
For instance,
the $L^2$-$m$-canonical forms
are just those $m$-canonical forms
on $X\backslash (D\cup\Sing X)$
that are locally $L^2$ at the
points of $D\cup \Sing X$,
see \cite{Sa}.  The sheaf
$\cL^{2,1}_X:=\cL^{2,1}_{X,0}$
was previously studied in
the paper of Grauert-Riemenschneider
\cite{GrRi} (it is called
there the canonical sheaf of $X$).

A motivation to study these
sheaves is that they allow simple
characterizations for when
a singularity is rational,
log terminal or log canonical.
As an example, by a result of Kempf \cite{KKMSD} a
normal complex algebraic variety $X$ has rational
singularities if and only if it is Cohen-Macaulay and
satisfies
$\cL^{2,1}_X\cong\cO_X(K_X)$.
In analogy with this result we will show:
if
$(X,x)$ is a normal complex singularity and $D$ is a reduced
Weil divisor with $K_X+D$ being
$\Q$-Cartier, then $(X,D)$ has
a log canonical singularity at $x$ if and
only if the stalks $(\cL^m_{X,D})_x$
and $\cO_X(m(K_X+D))_x$ are equal for all $m\ge 1$.
In the case when $D=0$ a similar
characterization holds for
log terminal singularities, see
\ref{Kempf}.

The main application of these
log canonical sheaves is to the case of
affine varieties $V=\Spec A$
for which the coordinate ring
$A=\bigoplus_{i\ge 0}A_i$ is
a non-negatively graded
$\C$-algebra.  Let $D\subseteq V$
be a reduced Weil divisor and assume
for simplicity that $D$ is
given by an equation $P=0$, where $P$ is
homogeneous of degree $d$.
The modules of sections
$$
L_{A,D}^{m}:=H^0(V,\cL_{V,D}^m)
\quad\mbox{and}\quad
L_{A,D}^{2,m}:=H^0(V,\cL_{V,D}^{2,m})
$$
as well as the reflexive hull
$\omega_A^{[m]}$ of the $m$-th tensor
power of the dualizing module $\omega_A$
then carry natural gradings.
We will show: {\em if
$V^*:= V\backslash V(A_+)$ is smooth and
$D\cap V^*$ is an SNC divisor then }
$$
L_{A,D}^{m}=(\omega_A^{[m]})_{\ge md}
\quad\mbox{ and }\quad
L_{A,D}^{2,m}=(\omega_A^{[m]})_{> (m-1)d}\,;
$$
see \ref{comprem} for a more general
statement.
This implies for instance that $(V,D)$
has log canonical singularities
if and only if $\omega_A^{[m]}$ has
no elements
of degree $< md$. Similar characterizations hold
for the properties
log terminal and rational, see \ref{cor logcan}.
For the latter case
this is a result of \cite{Fl} and \cite{KeWa}.

As another application we obtain
formulas for the
plurigenera of  the singularities of $V$,
where $V$ is as above.
For instance, the plurigenus $\delta_m$
introduced in
\cite{KiWa} is given by
$$
\delta_m (X,p)=\dim_{\C}\,(\omega_A^{[m]})_{\le
0},
$$
see \ref{del}. Applying this to complete
intersections we recover a result of
Morales \cite{Mo}.  A final
application concerns the logarithmic
plurigenera and the Kodaira
dimension of
$V\backslash D$. With the assumptions on $(V,D)$
as  above, assume moreover that $A_0\cong \C$
so that the $\C^*$-action
corresponding to the grading is good.
Then the logarithmic plurigenus
$\bar p_m(V\backslash D)$ is given by
the dimension of
$\omega^{[m]}_{md}$, see \ref{pluqua}.

The paper is organized as follows.
In Section 1 we introduce the sheaves of
log canonical forms and study
their basic properties. In particular
we show the characterization of
log canonical and log terminal singularities in
terms of  log canonical forms mentioned above.

Sections 2.2--2.5 contain the applications to
affine varieties with
$\C^*$-action as described above
In Section 2.1 we provide
some material concerning
equivariant completions
and weighted blowups of quasihomogeneous
affine varieties.

In this paper we work in the category
of complex spaces and
varieties over $\C$. However,
the principal results remain valid
for algebraic varieties over any
field of characteristic zero.

\section{Logarithmic $m$-canonical
forms on singular spaces}
\subsection{Logarithmic $m$-canonical
forms on smooth varieties}

\bnota
Let $X$ be a complex manifold
and $D\subseteq X$ a divisor
with simple normal crossings
(SNC in brief). Consider the
sheaf of logarithmic $m$-canonical forms
$$\cL_{X,D}^m:=\cO_X(mK_X+mD)\,$$
and the sheaf of logarithmic
$L^{2}$-$m$-canonical forms
$$\cL_{X,D}^{2,m}:=\cO_X(mK_X+(m-1)D)
=\cL_{X,D}^m(-D)\,.$$
If $x_1,\dots,x_n$ are local
coordinates around a point, say,
$p\in X$ with $D=\{x_1\cdot ... \cdot x_k=0\}$,
then near $p$ the
$\cO_X$-module $\cL_{X,D}^m$ is generated
by
$$
\omega_m=\left( \frac{dx_1}{x_1}
\wedge\dots\wedge
\frac{dx_k}{x_k}\wedge
dx_{k+1}\wedge\dots\wedge
dx_n\right)^{\otimes m}\,,
$$
whereas the sheaf $\cL_{X,D}^{2,m}$
is generated by $x_1\cdot ...
\cdot x_k\cdot \omega_m$.
As explained in \cite[Thm. 2.1]{Sa},
the forms in $H^0(X,\cL_{X,D}^{2,m})$
are just the meromorphic
$m-$canonical forms on $X\backslash D$
which locally at the
points of $D$ belong to $L^{2/m}$.

If $D=0$ is the zero divisor, we write $\cL_{X}^{m}$ and
$\cL_{X}^{2,m}$ instead of $\cL_{X,D}^{m}$,
resp.\ $\cL_{X,D}^{2,m}$.
\enota

\bsit\label{inj} Let $\pi:Y\to X$
be a morphism of complex manifolds of the same
dimension such that
$\pi^{-1}(D)_{\rm red}$ is contained
in an SNC divisor $E$ on
$Y$. Pulling back differential forms
induces natural
homomorphisms
$$\pi^*_m: \cL_{X,D}^m\to\cL_{Y,E}^m\qquad
{\rm and}\qquad
\pi^*_{2,m}: \cL_{X,D}^{2,m}\to\cL_{Y,E}^{2,m}\,
$$
(see \cite[\S 11.1.c]{Ii} for the case
of $\pi^*_m$; the case of
$\pi^*_{2,m}$ is similarly).
They are injective, if $\pi$ is dominant
\cite[Prop. 11.2]{Ii}. \esit

For later purposes we need the
following simple
observation.

\blem\label{sheaves} Let $\pi:Y\to X$
be a proper surjective
morphism of complex manifolds
of the same dimension, and
let $D\subseteq X$ be an SNC divisor
such that
$E=\pi^{-1}(D)_{\rm red}$ is also an
SNC divisor.  For a section $\eta\in
H^0(X\backslash D,\,\omega_X^{\otimes m})$
the following hold:
$$ \eta\in H^0(X,\,\cL_{X,D}^m)
\Leftrightarrow \pi^*\eta\in
H^0(Y,\,\cL_{Y,E}^m)$$
and
$$ \eta\in H^0(X,\,\cL_{X,D}^{2,m})
\Leftrightarrow \pi^*\eta\in
H^0(Y,\,\cL_{Y,E}^{2,m})\,.$$\elem

\proof The implication `$\Rightarrow$'
was already observed
before.  To show `$\Leftarrow$',
assume that $\pi^*\eta\in
H^0(Y,\,\cL_{Y,E}^m)$. Clearly, $\eta$
is a holomorphic section
of $\cL_{X,D}^m$ if it is locally
holomorphic outside an analytic
subset of $X$ of codimension 2.
Hence we may assume
that $\pi$ is finite and that
there are local coordinates
$x_1,\dots,x_n$ on $X$, resp.\
$y_1,\dots,y_n$ on $Y$ such that $D$
is given locally by $x_1=0$
and $\pi$ is given locally by $$
\pi:(y_1,y_2,\dots,y_n)\longmapsto
(y_1^k,y_2,\dots,y_n)
$$
for some $k\in \N$. The differential
form $\omega=\frac{dx_1}{x_1}
\wedge dx_2\wedge\dots\wedge  dx_n$
is then a local generator of the
invertible sheaf $\cO_X(K_X+D)$
and so $\omega^{\otimes m}$
generates $\cL_{X,D}^m$. Moreover,
$$\pi^*(\omega^{\otimes m})
=(k\frac{dy_1}{y_1}\wedge
dy_2\wedge\dots\wedge dy_n)^{\otimes m}\,$$
locally
generates $\cL_{Y,E}^m$, and
$\pi^*(x_1^a\omega^{\otimes m})
=y_1^{ak}\pi^*(\omega^{\otimes m}) $
is a local section of $\cL_{Y,E}^{2,m}$
if and only if $a\ge 1$.
This easily implies both statements
of the lemma.
\qed\medskip

Applying \ref{sheaves}
to local sections yields
the following corollary
(cf. \cite[Thm. 11.1]{Ii},
\cite[Thm. 1.1]{Sa}).

\bcor\label{IS} For a
bimeromorphic proper morphism $\pi:Y\to X$
of complex manifolds
and $D,\,E$ as above we have
$$
\pi_*(\cL_{Y,E}^m)=
\cL_{X,D}^m\qquad {\rm and}\qquad
\pi_*(\cL_{Y,E}^{2,m})=\cL_{X,D}^{2,m}\,.
$$
\ecor

\medskip

\subsection{Logarithmic $m$-canonical
forms on singular varieties}
In virtue of \ref{IS} the definition
of the sheaves $\cL_{X,D}^m$ and
$\cL_{X,D}^{2,m}$ can be
extended as follows.

\bdefi\label{ressheaf} Consider
a normal complex space $X$ and a
closed analytic subset $D\subseteq X$.
Let $\sigma:X'\to X$ be a
resolution of singularities such that
$D':=\sigma^{-1}(D\cup
\Sing X)_{\rm red}$ is an SNC divisor.
We call
$$\cL_{X,D}^m:=\sigma_*(\cL_{X',D'}^m)
\quad\mbox{and}\quad
\cL_{X,D}^{2,m}:=\sigma_*(\cL_{X',D'}^{2,m})\,$$
the sheaf of {\em logarithmic
$L^m$-canonical forms}, resp.\
{\em logarithmic $L^{2,m}$-canonical forms} on $X$.
\edefi

Because of \ref{IS} and by standard arguments,
this is independent
of the choice of resolution of singularities.
As before, if
$D=\emptyset$
then we write in brief $\cL_{X}^{m}$
and $\cL_{X}^{2,m}$
instead of
$\cL_{X,D}^{m}$ and $\cL_{X,D}^{2,m}$,
respectively.

\brem \label{1.6}
Clearly, $D\subseteq D_1$ implies that
$$\cL_{X,D}^{m}\subseteq
\cL_{X,{D_1}}^{m}\qquad {\rm and}\qquad
\cL_{X,D}^{2,m}\subseteq
\cL_{X,D_1}^{2,m}\,.$$
\erem

In most of our considerations
$D$ will be a Weil divisor. This is
justified by the following lemma.

\blem\label{divisorial part}
With $X$, $D$ as in
(\ref{ressheaf}) let
${\rm div} D$ denote the (reduced) union of all
divisorial components of $D$. Then the following hold.
\begin{itemize}
\item[ (a)] $\cL_{X,D}^m=\cL_{X,{\rm div} D}^m$ and
$\cL_{X,D}^{2,m}=\cL_{X,{\rm div} D}^{2,m}\,.$

\item[(b)] There are natural inclusions
$$
\cL_{X,D}^{m}\subseteq
\cO_X\,(m(K_X+{\rm div} D))
\quad\mbox{and}\quad
\cL_{X,D}^{2,m}\subseteq
\cO_X\,(mK_X+(m-1){\rm div} D)
$$
with equality outside the set  $\Sing X\cup \Sing {\rm div}
D$.
\end{itemize}
\elem

\proof
(a) We confine ourselves to
the proof of the first equality,
the proof of the second one being similarly.
By \ref{1.6} we have $\cL_{X,{\rm
div} D}^m\subseteq\cL_{X,D}^m$. To show the converse
inclusion, let $\pi:X'\to X$ be
a resolution of singularities such that
$\Delta':=\pi^{-1}({\rm
div} D\cup\Sing X)_\red$ and
$D':=\pi^{-1}(D\cup\Sing X)_\red$ are SNC
divisors on $X'$. Let $\eta$ be a local section of
$\cL_{X,D}^m$ in a
neighbourhood, say, $U$ of a point
$p\in X$.
Its restriction to
$U\backslash(D\cup\Sing
X)$ is a holomorphic section of
$\omega_X^{\otimes m}$ and so
extends holomorphically to an $m$-form
on $U\backslash({\rm
div} D\cup\Sing X)$. Thus $\pi^*(\eta)$ has no poles along
the components of $D'\backslash \Delta'$, so is a
section of $\cO_{X'}(m(K_{X'}+\Delta'))$, whence $\eta$ is a
section of $\cL_{X,{\rm div}D}^m$, as required.

In order to show (b), in view of (a), we may assume that $D$
is a reduced divisor. The sheaves $\cL_{X,D}^{m}$
and
$\cO_X\,(m(K_X+D))$
are then equal outside the
set $B:=\Sing X\cup\Sing D$.
As the latter sheaf is reflexive and  $\codim B\ge 2$ this
implies that
$$
\cL_{X,D}^{m}\subseteq
\cO_X\,(m(K_X+{\rm div} D)).
$$
The same argument also gives  the second inclusion in (b).
\qed\medskip

\brem\label{rem gen} More generally
one can introduce $m$- and
$\cL^{2,m}$-canonical forms for any
effective $\Q$-divisor
$D=\sum a_i D_i$ with $0<a_i\le 1$
\footnote{In \cite{FA} such a
divisor is called a {\it subboundary}.}.
Again, the definition is
given in two steps. If $X$ is a manifold
and $D_\red:=\sum D_i$
is an SNC divisor, we set
$$
\cL_{X,D}^m:=\cO_X( mK_X+\lfloor mD\rfloor )
\quad \mbox{and}
\quad \cL_{X,D}^{2,m}:=\cO_X( mK_X+\lfloor
(m-1)D\rfloor ).
$$
As $a_i\le 1$, with the same arguments as
above the construction is
functorial under generically finite maps
of manifolds of the same
dimension.

In the case of a normal variety $X$ we
choose a resolution
of singularities $\sigma:X'\to X$
such that $\sigma^{-1}(D\cup\Sing X)_\red$ is a simple
normal crossing divisor. Let $D'$
denote the divisor $D^\pr+
\sum_iE_i$, where $E_i$ are the exceptional
divisors and $D^\pr$ denotes
the proper transform of $D$. Now one can
introduce as in \ref{ressheaf}
the sheaves
$$
\cL_{X,D}^m:=\sigma_*(\cL_{X',D'}^m)
\quad\mbox{and}\quad
\cL_{X,D}^{2,m}:=\sigma_*(\cL_{X',D'}^{2,m}).
$$
As before one can show that this
definition does not depend on the
choice of the resolution.
\erem

The following proposition indicates
certain functorial properties
of $L^m$- and $L^{2,m}$-canonical forms
(cf. \cite[Prop. 11.3]{Ii}).

\bprop\label{dirim}
Let $\pi:Y\to X$ be a generically finite
morphism of normal connected complex spaces
of the same dimension. Let
$D\subseteq X$ be an analytic subset, and
assume that $E\subseteq Y$ is
an analytic subset with  ${\rm div}E={\rm
div}\pi^{-1}(D\cup\Sing X)$. Then the following hold.

\begin{enumerate}
\item[(a)] There are natural injections
$$
\cL_{X,D}^m\to\pi_*(\cL_{Y,E}^m)\qquad
      {\rm and}\qquad \cL_{X,D}^{2,m}\to
\pi_*(\cL_{Y,E}^{2,m})\,.$$

\item[(b)] If moreover $\pi$ is proper,
then for a form
$\eta\in H^0(X\backslash D,\,\cO_X(mK_X))$
we have
$$
\begin{array}{rcl}
\eta\in H^0(X,\,\cL_{X,D}^m)&\Leftrightarrow&
\pi^*\eta\in
H^0(Y,\,\cL_{Y,E}^m)
\quad\hbox{and}\quad\\[2pt]
\eta\in H^0(X,\,\cL_{X,D}^{2,m})
&\Leftrightarrow& \pi^*\eta\in
H^0(Y,\,\cL_{Y,E}^{2,m})\,.
\end{array}
$$

\item[(c)] If $\pi$ is proper and
birational,
then
$$\cL_{X,D}^m=\pi_*(\cL_{Y,E}^m)\qquad
      {\rm and}\qquad \cL_{X,D}^{2,m}
=\pi_*(\cL_{Y,E}^{2,m})\,.
$$\end{enumerate}
\eprop

\proof
Because of \ref{divisorial part} we may assume that
$E=\pi^{-1}(D\cup \Sing X)$.  As before we consider
only the case of
$\cL^m$-forms, the
case
of $\cL^{2,m}$-forms being similarly.
Consider resolutions of singularities
$Y'\to Y$ and
$X'\to X$ that fit into a diagram
\begin{diagram}[s=7mm]
Y'&\rTo^{\pi'} & X'\\
\dTo<q & & \dTo>p \\
Y & \rTo^\pi & X
\end{diagram}
and such that
$$
E':=q^{-1}(E \cup {\rm Sing}\,Y)_\red\quad
\mbox{and}\quad
D':=p^{-1}(D\cup   {\rm Sing}\,X)_\red
$$
are SNC divisors in $Y'$ resp.\ $X'$.
As $\pi^{\prime-1}(D')$ is
contained in $E'$ the morphism $\pi'$
induces an injection
$\cL_{X',D'}^m\to\pi'_*(\cL_{Y',E'}^m)$
(see \ref{inj}).
Applying $p_*$
gives the desired injection in (a).

(c) is an immediate consequence of (b).
To deduce (b), note
first that $\pi(\Sing Y)\subseteq X$ is
a closed analytic subset of
codimension at least 2. By \ref{divisorial
part} the sheaves
$\cL_{X,D}^m$ and $\cL_{X,D_1}^m$ are equal,
where
$D_1:=D\cup\pi(\Sing Y)$. Moreover $\Sing Y$
is contained in
$E_1:=\pi^{-1}(D_1)$ and $\cL_{Y,E}^m
\subseteq \cL_{Y,E_1}^m$.
Hence it is sufficient to prove (b) for
$D_1$ instead of $D$.
In other words, we may assume
that $E$ contains $\Sing Y$.

Let now $q:Y'\to Y$ and  $p:X'\to X$
be resolutions of
singularities as
above. By \ref{sheaves}
$$
p^*\eta\in H^0(X',\,\cL_{X',D'}^m)
\Leftrightarrow
(\pi'p)^*\eta\in
H^0(Y',\,\cL_{Y',E'}^m).
$$
As by definition $p_*(\cL_{X',D'}^m)
=\cL_{X,D}^m$ and
$q_*(\cL_{Y',E'}^m)
=\cL_{Y,E}^m$\, , (b) follows.
\qed\medskip

For our purposes it is useful
to introduce
certain sheaves that are invariants of the
singularities.

\bnota \label{rn}  Let $X$, $D$ be as
in \ref{ressheaf}. Because of \ref{divisorial part} (b) we
may form the quotient sheaves
$$\gL^m_{X,D}:=\cO_X\,(m(K_X+\,{\rm div}
D))\bigr/\cL_{X,D}^{m}$$
and $$\gD^{m}_{X,D}:=
\cO_X\,(mK_X+(m-1){\rm div} D)
\bigr/\cL_{X,D}^{2,m}\,.$$
Note that by \ref{divisorial part} (a)
$\gL^m_{X,D}=\gL^m_{X,{\rm div}D}$ and
$\gD^m_{X,D}=\gD^m_{X,{\rm div}D}$. Moreover, by
\ref{divisorial part} (b) these sheaves are concentrated on
$\Sing X\cup \Sing {\rm div}D$. As before, in the case
$D=\emptyset$
we write in brief
$\gL^m_X$ and $\gD^m_X$ instead of $\gL^m_{X,D}$,
resp.\ $\gD^m_{X,D}$.
\enota

Later on we will need the following fact.

\blem\label{injective}
For analytic subsets $D_1\subseteq D_2$ of $X$ the
natural maps
$$
\gL_{X,D_1}^m\lto \gL_{X,D_2}^m
\quad\mbox{and}\quad
\gD_{X,D_1}^m\lto \gD_{X,D_2}^m
$$
are injective.
\elem

\begin{proof}
We restrict to the proof of the first
inclusion the other one being similarly.
By \ref{divisorial part} (a) we may assume that $D_1$ and
$D_2$ are reduced divisors.  We need to show that
$$
\cL_{X,D_1}^m=\cL_{X,D_2}^m\cap \cO_X(m(K_X+D_1)).
$$
The inclusion `$\subseteq$' follows from \ref{1.6} and
\ref{divisorial part} (b). To show the converse inclusion,
let
$\pi:X'\to X$ be a resolution of singularities such that
$D'_i:=\pi^{-1}(D_i\cup \Sing X)_\red$, $i=1,2$, are SNC
divisors. If $\eta$ is a section of $\cL_{X,D_2}^m\cap
\cO_X(m(K_X+D_1))$ defined over some open subset $U$ of $X$,
then $\pi^*(\eta)$ is a form in $\omega_{X'}^{\otimes m}$
that has at most logarithmic poles along the irreducible
components of $D'_2$ and is holomorphic along the components
of the proper transform of $D_2- D_1$. Hence
$\pi^*(\eta)$ is a section of
$\cO_{X'}(m(K_{X'}+D_1'))$ so that $\eta\in
H^0(U,\cL^m_{X,D_1})$, as required.
\end{proof}

In the next proposition we
describe the behavior of $L^m$- and
$L^{2,m}$-canonical forms and
the sheaves $\gL^m$ and $\gD^m$
under taking Cartesian products.

\bprop\label{products}
Let $X_1$, $X_2$ be normal complex spaces, let
$D_i\subseteq X_i$,
$i=1,2$, be closed analytic subsets and
let $D$ denote the closed
analytic subset $X_1\times
D_2 \cup D_1\times X_2$ of the product
$X:=X_1\times X_2$.
With $p_i:X\to X_i$ being the canonical
projection
($i=1,2$) the following
hold.
\smallskip

(a) The sheaves $\cL^m$, $\cL^{2,m}$,
$\gL^m$ and
$\gD^m$ ($m\ge 1$) are compatible with taking products, i.e.\ with
$\cH$ any one of these sheaves we have
$$
\cH_{X,D}\cong p_1^*(\cH_{X_1,D_1})
\otimes
p_2^*(\cH_{X_2,D_2}).
$$

(b) If $D_1$ and $D_2$ are divisors then for all $k,m\in \Z$
$$
\cO_X(mK_X+kD)\cong p_1^*(\cO_{X_1}(mK_{X_1}
+kD_1))
\otimes
p_2^*(\cO_{X_2}(mK_{X_2}+kD_2))\,.
$$

(c)  If $D_1$, $D_2$, are divisors
then $D$ is
$\Q$-Cartier if and only
if $D_1$ and $D_2$ are $\Q$-Cartier.
\eprop

\proof
(b) is obvious on the regular part of $X$. The sheaves
on both sides of (b) are reflexive and so
are determined by their restrictions to the regular parts,
whence (b) follows.

In order to deduce (a),
let $\pi_i:X_i'\to X_i$ ($i=1,2$)
be resolutions of
singularities such
that $D_i':=\pi_i^{-1}(D_i\cup \Sing X_i)$
are SNC divisors in
$X_i'$. The
product $X':=X_1'\times X_2'$ then provides
a resolution of
singularities
$\pi:X'\to X$ such that $D':=\pi^{-1}(D\cup
\Sing
X)=X_1'\times D_2' \cup D'_1\times X'_2$ is
an SNC divisor.
Let
$q_i:X'\to X_i'$ denote the canonical
projection ($i=1,2$).
By (b)
$$
\cL_{X',D'}^m\cong q_1^*(\cL_{X_1',D_1'}^m)
\otimes
q_2^*(\cL_{X_2',D_2'}^m)\quad\mbox{and}\quad
\cL_{X',D'}^{2,m}\cong q_1^*(\cL_{X_1',D_1'}^{2,m})
\otimes
q_2^*(\cL_{X_2',D_2'}^{2,m}).
$$
Applying $\pi_*$ and using the K\"unneth formula
gives (a) for the cases  $\cH=\cL^m$ and $\cH=\cL^{2,m}$.
Using (b) also the remaining two cases follow by taking
quotients.

(c) is an immediate consequence of the fact
that for $n\in \N$
the divisor
$nD$ is Cartier if and only if $nD_1$ and
$nD_2$ are
Cartier.
\qed\medskip

For our applications to quasihomogeneous
singularities
it is important
to study the behaviour of $L^m$- and
$L^{2,m}$-canonical forms
under finite group actions.

\bprop\label{fingr} Let $G$ be a finite
group acting
on a normal complex space $Y$, and let
$\pi:Y\to X:=Y/G$ be the canonical
morphism onto the orbit space. Let
$D\subseteq X$ be an
analytic subset
and assume that
$\pi$ is unramified in codimension
one outside
$E:=\pi^{-1}(D)_{\rm
red}$. Then the following hold.

\begin{enumerate}
\item[(a)]
$\,\,\cL_{X,D}^{m}=\pi_*(\cL_{Y,E}^{m})^G$
and
$\cL_{X,D}^{2,m}=\pi_*(\cL_{Y,E}^{2,m})^G\,.$
\item[(b)] $\gL^m_{X,D} = (\pi_*\gL^m_{Y,E})^G$
and
$\gD^m_{X,D} = (\pi_*\gD^m_{Y,E})^G$.
\item[(c)] If $D$ is a divisor, then $\pi^*(D)$
is $\Q$-Cartier if and only if $D$ itself is
$\Q$-Cartier.
\end{enumerate}
\eprop

\proof
In order to show (a) we may assume
that $D$, and then also
$E$, are Weil divisors, see
\ref{divisorial part}. With
$E':=\pi^{-1}(D\cup \Sing X)$,
\ref{dirim} (b) implies that
$\cL_{X,D}^{m}=\pi_*(\cL_{Y,E'}^{m})^G$
and
$\cL_{X,D}^{2,m}=\pi_*(\cL_{Y,E'}^{2,m})^G.$
As $E$ and $E'$
are equal in codimension 1, (a) follows from
\ref{divisorial part}.

(a) implies in particular that
\begin{eqnarray}
\label{e1} \cO_X\,(m(K_X+D))
&=&\pi_*(\cO_Y\,(m(K_Y+E)))^G\\
\label{e2} \cO_X\,(mK_X+(m-1)D)
&=&\pi_*(\cO_Y\,(mK_Y+(m-1)E))^G\,.
\end{eqnarray}
Indeed, the involved sheaves are reflexive and so it
is sufficient to verify the equalities on the part
where $X$ and $Y$ are both smooth. Now (b) follows
from (a) in view of (\ref{e1}) and (\ref{e2}). The
statement of (c) is well  known
(see e.g., \cite{St, Fo}).
\qed\medskip

Recall the following notion.

\bdefi A morphism of complex spaces
$\pi:Y\to X$ is called {\em
non-degenerate} if for every point
$y\in Y$ we have $\dim_yY=\dim_y
\pi^{-1}(\pi(y))+\dim_{\pi(y)}X$.
\edefi

For instance, every finite
surjective morphism is non-degenerate.
The following proposition will be
used in the sequel.

\bprop\label{inclusion}
Let $\pi:(Y,y)\to (X,x)$ be a non-degenerate
morphism of
normal complex space germs, and let
$D\subseteq X$ be a reduced Weil
divisor with preimage $E:=\pi^{-1}(D)_{\rm red}$.
Then there are (in
general non-canonical) injections of
$\cO_{X,x}$-modules
\be\label{gL} (\gL^m_{X,D})_x \hookrightarrow
(\gL^m_{Y,E})_y\ee
and
\be\label{gD} (\gD^m_{X,D})_x \hookrightarrow
(\gD^m_{Y,E})_y\,.\ee
\eprop

\proof
We restrict ourselves to the proof (\ref{gD}), the other one
being similarly. First we treat the special case
that $\pi$ is finite. Pulling back differential forms
induces an injective map
\be\label{eq 5}
\pi^*:\cO_X(mK_X+(m-1)D)\hto
\cO_Y(mK_Y+(m-1)E).
\ee
The analytic sets
$E$ and $\pi^{-1}(D\cup \Sing X)_\red$
are equal in codimension 1 and so  by \ref{dirim} (b)
for a local section, say, $\eta$ of $
H^0(U,\cO_X(mK_X+(m-1)D))$ over some open subset $U\subseteq
X$
$$
\eta\in H^0(U,\cL_{X,D}^{2,m})
\quad\mbox{if and only if}\quad
\pi^*(\eta)\in H^0(\pi^{-1}(U),\cL_{Y,E}^{2,m}).
$$
Thus (\ref{eq 5}) induces an injective map as in (\ref{gD}).

In the general case we can find functions
$f_1, \ldots,f_d\in\cO_{Y,y}$
vanishing at $y$, where $d:=\dim_y \pi^{-1}(x)$,
such that
$f:=(f_1,\ldots,f_d)$ restricts to a finite
map of germs
$(\pi^{-1}(x),y)\to (\C^d,0)$. Thus, letting
$Z:=X \times\C^d$ and $z:=(\pi(y),0)$,
we can factorize
$\pi$ into two maps
$$
\pi: (Y,y)\stackrel{\pi\times f}{\lto}(Z,z)
\stackrel{pr_1}{\lto} (X,x),
$$
where $\pi\times f$ is finite.
With $D_Z:=D\times \C^d$ we have
(see \ref{products} (a))
\be\label{gD1}
(\gD^m_{Z,D_Z})_z\cong (\gD^m_{X,D})_x
\otimes_{\cO_{X,x}}\cO_{Z,z},
\ee
Applying
the first part of the proof, (\ref{gD}) follows.
\qed\medskip

\subsection{Logarithmic $m$-canonical forms versus
log canonical, log terminal and rational
singularities}
Recall the following notions.

\bdefi\label{lc} Let $X$ be a normal
complex space and let
$D$ be a reduced effective Weil divisor
on $X$ such that
the divisor $K_X+D$ is $\Q$-Cartier.
Let $\sigma:X'\to X$ be a
resolution of singularities such
that $D'=\sigma^{-1}(D\cup \Sing
X)_{\rm red}$ is an SNC divisor.
Write
$$
\sigma^*(K_X+D)=(K_{X'}+D')-\sum_i
a_iE_i\qquad
{\rm with}\quad a_i\in\Q\quad\forall i,
$$
where the summation is taken over the
set of all divisorial
irreducible components $E_i$ of the
exceptional set of the blowup
$\sigma$. (The number $a_i$ is the
so called {\em log
discrepancy}\/ of $E_i$ \cite[2.5.3]{FA}.)
One says\footnote{See
e.g. \cite[Def. 2.34]{KoMo},
\cite[Def. 3.5]{Ko}.}
that the pair
$(X,\,D)$ has

\begin{enumerate}
\item[1.] {\it log canonical singularities}
if $a_i\ge 0\quad\forall i$.
\item[2.] {\it log terminal
singularities}\footnote{In \cite{KoMo}
these singularities are
called {\em purely log terminal.}}
if $a_i>0\quad\forall i$ for {\em every}\/
resolution $\sigma:X'\to X$.
\end{enumerate}

Following \cite{Ko} we will simply say that
$(X,\,D)$ is log canonical (lc, for short)
resp., log terminal
(lt, for short). In the case when $D$ is
the zero divisor
one says that $X$
(instead of $(X,0)$) has log canonical
resp.\ log terminal
singularities.

\begin{enumerate}
\item[3.] (Artin \cite{Ar})
$X$ is said to have {\em rational
singularities} if for a
resolution of singularities $\sigma:X'\to X$
the higher direct
image sheaves $R^i \sigma_*(\cO_X)$, $i\ge 1$,
vanish.
\end{enumerate}
\edefi

In the next proposition we recall Kempf's
characterization
of rational singularities in terms of
canonical forms. Moreover,
we characterize log canonical and log terminal
singularities in terms of logarithmic
$m$-canonical forms.

\bprop\label{Kempf}
\begin{enumerate}
\item[(a)] {\em (Kempf \cite[Prop.\ on
p.\ 50]{KKMSD})}
$X$ has rational singularities if
and only if $X$
is Cohen-Macaulay and
$\cL^{2,1}_{X}\cong \cO_X(K_X)$ i.e.,
$\gD^1_X=0$.
\end{enumerate}

Assume further that $K_X+D$ is a\/
$\Q$-Cartier
divisor. Then

\begin{enumerate}
\item[(b)] $(X,\,D)$
is lc if and only if
$\cL^m_{X,D}=\cO_X(m(K_X+D))$ for all
$m\ge 1$ or, equivalently, $\gL^m_{X,D}=0$.
\item[(c)] $X$ is lt if and only if
$\cL^{2,m}_{X,D}=\cO_X(mK_X+(m-1)D)$ for all
$m\ge 1$ or, equivalently, $\gD^m_{X,D}=0$.
\end{enumerate}
\eprop

\proof (a) The original result of
Kempf (formulated only for
algebraic singularities) generalizes
to the complex analytic
setting in view of the fact that the
Grauert-Riemenschneider
vanishing theorem \cite{GrRi} also
holds for non-algebraic
singularities. The latter follows
from a more general result of
Moriwaki \cite[Thm. 3.2]{Mor}.

\no (b) By assumption $k(K_X+D)$ is
a Cartier divisor for some
$k\in\N$. Fix a resolution of
singularities $\sigma:X'\to X$ such
that $D':=\sigma^{-1}(D\cup \Sing X)_{\rm red}$
is an
SNC divisor. By \ref{lc}, $(X,\,D)$ is lc
if and only if
\be \label{logcan}
\cO_{X'}\,(\sigma^*(k(K_X+D)))\subseteq
\cO_{X'}\,(k(K_{X'}+D'))\,. \ee Thus for
a local section $\omega$
of $\cO_X(m(K_X+D))$ the tensor power
$\omega^{\otimes k}$ pulls
back to a section $\sigma^*\omega^{\otimes k}$ of
$\cO_{X'}(km(K_{X'}+D'))$, whence
$\sigma^*(\omega)$ is a section
of $\cO_{X'}(m(K_{X'}+D'))$ and so by definition
(see \ref{ressheaf})
$$\cO_X(m(K_X+D))\subseteq
\cL_{X,D}^{m}\quad\mbox{for all }m\ge 1\,.$$
Since $\cL_{X,D}^{m}\subseteq
\cO_X(m(K_X+D))$ (see \ref{divisorial part} (b)), the latter
means that
\be\label{star}
\gL_{X,D}^{m}=0 \qquad \mbox{for all } m\ge 1\,.
\ee
Conversely, if $\gL_{X,D}^{k}=0$  then
(\ref{logcan}) holds, proving (b).
The proof of (c) is similarly
and left to the reader.\qed\medskip

For further purposes, it is convenient
to introduce the following definition.

\bdefi\label{L2} We will say that the pair
$(X,D)$ is $L^2$-{\em log terminal} ($L^2$-lt,
for short) if
$\gD^m_{X,D}=0\quad\forall m\ge 1$. \edefi

\bprop\label{l2t} Let $X$ be a normal
complex space and $D$ be a
divisor on $X$ such that $\Sing
X\subseteq D$. Assume that both
$K_X$ and $D$ are $\Q$-Cartier.
If $(X,D)$ is lc, then it is $L^2$-lt.
\eprop

\no\proof Let as above $\sigma:X'\to X$
be a resolution of
singularities such that $D':=\sigma^{-1}(D)_\red
=\sigma^{-1}(D\cup
\Sing X)_{\rm red}$ is an SNC divisor.
We need to show that
$\cO_X\,(mK_X+(m-1)D)\subseteq \cL_{X,D}^{2,m}$
for all $m\in\N$
or, equivalently, that for every local section
$\omega$ of
$\cO_X\,(mK_X+(m-1)D)\subseteq \cO_X\,(m(K_X+D))$
the form
$\sigma^*(\omega)$ extends to a section in
$\cL_{X',D'}^{m}(-D')$.
By assumption, (\ref{star}) holds and
therefore
$\sigma^*(\omega)$ gives a section in $\cL_{X',D'}^{m}$.
Choose
$k$ such that $kD$ is a Cartier divisor, so that
$kD$ is locally
given on $X$ by one equation, say $f=0$.
The form
$\omega^{\otimes k}$ is a section in
$\cO_X\,(kmK_X+k(m-1)D)=\cL_{X,D}^{km}(-kD)$,
and so it can be
written locally as $f\cdot\eta$ for
some section $\eta$ in
$\cL_{X,D}^{km}$. Hence $\sigma^*(\omega^{\otimes
k})=\sigma^*(f\eta)
=\sigma^*(f)\cdot\sigma^*(\eta)$ becomes a
section of $\cL_{X',D'}^{km}$ that
vanishes along $D'$. It
follows that $\pi^*(\omega)$ as a section of
the sheaf $\cL_{X',D'}^{m}$
vanishes along $D'$ as well, which
implies the assertion.\qed\medskip

\begin{rem}\label{rem gen2}
(1) Thus for pairs $(X,D)$ satisfying
the assumptions of
\ref{l2t} the
following inclusions hold:
$${\rm (lt)} \subseteq {\rm (lc)}
\subseteq {\rm (}L^2-{\rm lt)}\,,$$
whereas for $D=0$ by \ref{Kempf} (c)
we have:
$${\rm (lt)}={\rm (}L^2-{\rm lt)}
\subseteq {\rm (lc)}\,.$$
The latter equality is no longer
true for pairs $(X,D)$
if $D\ne 0$. The
simplest example is given by the
union $D$ of two smooth curves
on a smooth surface $X$ meeting
transversally. The pair $(X,D)$
is not lt although it is lc and
hence (by \ref{l2t}) it is $L^2$-lt.

(2) Using the definition of log
canonical forms given in
\ref{rem gen} one can also
characterize log canonical
singularities if
$D=\sum a_iD_i$ is a $\Q$-divisor with
$0<a_i\le 1$. Moreover, using
this characterization one can extend
the notion of log canonical
singularities without requiring that
$K_X+D$ is
$\Q$-Cartier.
\end{rem}

Combining \ref{dirim} (c) and \ref{Kempf}
gives the following corollary.

\bcor\label{log can}
Let $\pi:Y\to X$ be a proper surjective
bimeromorphic morphism of connected
normal complex spaces
and let $D\subseteq X$ be an analytic
subset.
Denote by $E$ the union of divisorial
components
of the analytic subset
$\pi^{-1}(D\cup {\rm Sing}\,X)_{\rm
red}$.

\begin{enumerate}
\item[(a)] If $(Y,\,E)$ is lc then
$$\cL_{X,D}^{m}=\pi_*(\cO_Y\,(m(K_Y+E)))\,.$$

\item[(b)] If $(Y,E)$ is $L^2$-lt
then
$$
\cL_{X,D}^{2,m}=\pi_*(\cO_Y\,
(mK_Y+(m-1)E))\,.$$
\end{enumerate}\ecor

   From \ref{fingr} and \ref{Kempf}
we obtain the following corollary.

\bcor \label{fingrcor} Let $\pi:Y\to X=Y/G$,
$D\subseteq X$ and $E\subseteq Y$
be as in \ref{fingr}.

\begin{enumerate}
\item[a)] If $(Y,E)$ is $L^2$-lt then
so is $(X,D)$.
\item[b)] If  $D$ is a reduced divisor and
$\pi^*(K_X+D)$ is
$\Q$-Cartier then the same is true
for the property  `lc'.
\end{enumerate} \ecor

The next corollary follows from
\ref{products} (a),(c)
using \ref{Kempf} (b).

\bcor Under the assumptions as in
\ref{products}
if $D_1,\,\,D_2$ are divisors then
$(X,D):=(X_1\times X_2, X_1\times D_2
\cup D_1\times X_2)$
is lc (resp., $L^2$-lt) if and only if
$(X_1,D_1)$ and $(X_2,D_2)$ are lc
(resp., $L^2$-lt).\ecor

The following corollary is an immediate
consequence
of \ref{inclusion}
and \ref{Kempf}. Part (a) was shown by
Bingener and Storch, whereas part (b)
generalizes a result of Ishii and Koll\'ar.

\bcor\label{cor kollar}
For a  non-degenerate surjective morphism of normal
complex spaces $\pi:Y\to X$ the following hold.
\begin{enumerate}
\item[(a)] \cite[5.7]{BiSt} If $Y$ has
rational singularities
and $X$ is Cohen-Macaulay, then
$X$  has also rational singularities.
\item[(b)] {\rm (cf. \cite[1.7.II]{Ish}, \cite[20.3.3]{FA})}
Let $D\subseteq X$
be a reduced Weil divisor with
preimage $E:=\pi^{-1}(D)_{\rm red}$.
If $(Y,E)$ is $L^2$-lt then so is $(X,D)$.
If $K_X+D$ is $\Q$-Cartier, then the
same holds for the property `lc'.
\end{enumerate}
\ecor

\subsection{Plurigenera of isolated
singularities}
Let $X$ be a normal complex space.
Recall the following notions and
facts \cite{KiWa,
Kn, Mo}.

\bdefi\label{pgr} {\rm
By \cite{KiWa}, the sections of
the sheaf $\cL_{X}^{2,m}$
over an open subset $U\subseteq X$
are just the sections in
$H^0(U_{\rm reg},\,\omega_X^{\otimes m})$
which are locally
$L^{2/m}$ on $U$.
If $(X,\,x)$ is an isolated singularity
then
$$
\delta_m\,(X,\,x) =
\dim_{\C}\,\left[\omega_{X,\,x}^{[m]}
/(\cL_{X}^{2,m})_x\right]\,
=\dim_{\C}\,(\gD^m_X)_x$$
is the {\it $m$-th $L^2$-plurigenus} as
defined in
\cite{KiWa}; here  $\omega_X^{[m]}:=\cO_X(mK_X)$.
Similarly, the  {\it $\lambda$-plurigenera}
\cite{Mo} are given by
$$
\lambda_m\,(X,\,x) =
\dim_{\C}\,\left[\omega_{X,\,x}^{[m]}/
(\cL_{X}^{m})_x\right]=
\dim_{\C}\,(\gL^m_X)_x\,.$$}\edefi

\bsit\label{kmp} Note that by Kempf's
criterion \ref{Kempf},
{\it an isolated  singularity $(X,x)$
is rational if and only if it is
Cohen-Macaulay and $\delta_1(X,x)=0$}.
Moreover, {\it if $(X,x)$ is
$\Q$-Gorenstein\footnote{Recall \cite[1.3-1.4]{Ish}
that $(X,x)$ is Gorenstein iff it is Cohen-Macaulay and
the canonical sheaf $\omega_X$ is invertible at $x$. It is called
$\Q$-{\it Gorenstein} if the sheaf $\omega_X^{[m]}$
is invertible at $x$ for some $m\in\N$.} then by
\ref{Kempf} $(X,x)$ is lt (resp., lc)
if and only if all plurigenera $\delta_m$ (resp.\
$\lambda_m$), $m\ge 1$, vanish.}\esit

In analogy with \ref{cor kollar}
we have the following result.

\bcor\label{cor plurignera}
If $\pi:(Y,y)\to (X,x)$ is a non-degenerate
surjective morphism of
germs of normal isolated singularities then
$$
\delta_m(X,x)\le \delta_m(Y,y)
\quad \mbox{and}\quad
\lambda_m(X,x)\le \lambda_m(Y,y)
\quad\forall m\ge 1\,.$$
Moreover, if $\dim_yY>\dim_xX$ then
$$
\delta_m(X,x)=\lambda_m(X,x)=0
\quad\forall m\ge 1\,.$$
\ecor

\proof
The first part follows from the inclusions
(\ref{gL}) and (\ref{gD})
in \ref{inclusion}. To
show the second assertion, consider a
factorization  $(Y,y)\to (Z,z)\to
(X,x)$  as in the proof of \ref{inclusion}.
Assuming that
$\delta_m(X,x)\neq 0$, by (\ref{gD1}) in the proof of {\it
loc.cit.}\ the module
$(\omega_Z^{[m]}/\cL_{Z}^{2,m})_z=(\gD^m_Z)_z$
has infinite dimension and is contained in
the finite
dimensional vector space
$(\omega_Y^{[m]}/\cL_{Y}^{2,m})_y=(\gD^m_Y)_y$.
Thus we get a contradiction, and so
$\delta_m(X,x)$
vanishes. As $\lambda_m\le \delta_m$,
the $\lambda$-plurigenera vanish as well.
\qed\medskip

\section{Logarithmic $m$-canonical forms on
quasihomogeneous varieties}
In this section we show how to compute
$L^m$- and $L^{2,m}$-canonical forms
on affine varieties with
$\C^*$-action, and we apply this to
characterize different types of singularities.

\subsection{Graded rings and associated
schemes}\label{gras}

For the convenience of the reader we recall
some facts about
projective schemes associated to graded
rings which will be useful
in the sequel (see \cite[sect. 2]{Fl}).

\bnota\label{a.1}
Let $K$ denote a field of characteristic 0
containing all
roots of unity. Recall that for a finitely
generated graded $K$-algebra
$A=\bigoplus_{\nu\ge 0} A_{\nu}$ the
associated projective scheme
$\Proj A$ is defined by the set of all
homogeneous
prime ideals
$\fp$ of $A$ with $A_+\not\subseteq \fp$,
where $A_+=\bigoplus_{\nu> 0} A_{\nu}$
is the augmentation ideal. The scheme $\Proj A$
is separated and
proper over $\Spec A_0$
\cite{EGA}. Furthermore,
$\Proj A$ is covered by the affine open subsets
$$
D_+(f)=D_+(fA):=\{\fp\in \Proj A :
f\notin \fp\}\cong \Spec A_{(f)}\,,
$$
where $f\in A_+$ is a homogeneous
element and $A_{(f)}:=(A_{f})_0$
denotes the degree zero part of the
localization $A_f$.
Denote also $V=\Spec A$ and
$V^*=V \backslash V(A_+)$
(where, as usual, $V(I)$ is the zero
set of an ideal $I$, whereas for
a homogeneous ideal $I\subseteq A$, $V_+(I)$
stands for its zero set in $\Proj  A$).
There is a natural surjective morphism
$V^*\to \Proj  A$.
\enota

The multiplicative group $\G_m$ of
the field $K$ acts
on $A$ via
$t.a=t^{\nu}a$ for $t\in K^*$ and
$a\in A_{\nu}$.
If $A=A_0[A_1]$ is
generated as an $A_0$-algebra by the elements
of degree
1 then $V^*\to \Proj A$ is a locally trivial
$\G_m$-bundle.
In general, we have the following well known
fact.

\blem\label{sl}
$\Proj  A\cong V^*/\G_m$.
\elem

\proof
In lack of a reference we provide the
simple argument: $V^*$ is covered by the
$\G_m$-invariant affine open subsets
$D(f):=\Spec A_f$, where
$f\in A_{d}$ with $d >0$.
As the ring of invariants
$(A_f)^{\G_m}$ is just $A_{(f)}=(A_f)_0$ we obtain
$D_+(f)=D(f)/{\G_m}$ and so the lemma follows.
\qed\medskip

\bsit To describe the situation more closely, for
$f\in A_{d}$ with $d>0$ denote $F=F(f)=A/(f-1)$
resp.,
$Y=Y(f)=\Spec F$
and consider  the homomorphism of graded rings
$$
\mu : A_f\to F[T,T^{-1}]
\quad\mbox{given\,\,\, by}\quad
a/f^k\longmapsto {\overline a}\cdot
T^{\deg a-kd}\,,
$$
where $F[T,T^{-1}]$ is graded via
$F[T,T^{-1}]_0=F$ and $\deg T=1$.
Clearly
$\mu$ is degree preserving; in particular
$\mu ((A_f)_0)\subseteq F$.
There is a commutative diagram
\begin{diagram}[w=7mm,h=8mm,midshaft]
A_{(f)}=&(A_f)_0 & \rTo &  F     \\
           &\dTo<i  &            &  \dTo>i\\
           &A_f     & \rTo^{\mu}       &  F[T,T^{-1}]
\end{diagram} where $i$ stands for the natural
inclusions.

The cyclic group
$\Z_{d}\cong\langle\zeta\rangle$ generated by a
primitive
$d$-th root of unity $\zeta$ acts (homogeneously)
on $F[T,T^{-1}]$;
namely for
$b={\overline a}T^k\in F[T,T^{-1}]$ with
$a\in A$ homogeneous we let
$\zeta.b=\zeta^{\deg a-k}b\,.$ This action
restricts to
$F=F[T,T^{-1}]_0$;
the next lemma describes the quotients. \esit

\blem\label{il} \cite[2.1-2.2]{Fl} (a)
$\mu$ provides isomorphisms
$A_f\cong F[T,T^{-1}]^{\Z_{d}}$ resp.,
$A_{(f)}=(A_f)_0\cong F^{\Z_{d}}$
onto the rings of invariants, and the horizontal
arrows
(i.e., the orbit maps)
in the induced commutative diagram
\begin{diagram}[w=7mm,h=8mm,midshaft]
\Spec F[T,T^{-1}] &= & Y\times \G_m
& \rTo^{\ul\mu}_{/\Z_d}& D(f)& =&
\Spec A_f     \\
                       &&\dTo<{{\rm pr}_1}   &
&\dTo>{/\G_m}\\
    \Spec F           &=&  Y
& \rTo_{/\Z_d}
& D_+(f)& =&
\Spec A_{(f)}
\end{diagram} are cyclic coverings; moreover,
${\ul\mu}$ is an \'etale covering.

\smallskip\no  (b) Furthermore, $\mu$
induces isomorphisms
$$
(A_f)_{\ge 0}\cong F[T]^{\Z_{d}}
\quad\mbox{and}\quad
(A_f)_{\le 0}\cong F[T^{-1}]^{\Z_{d}}\,,
$$
and so it provides
ramified cyclic coverings
$$Y\times
\A^1_K \stackrel{/\Z_d}{\longrightarrow}
\Spec\, (A_f)_{\ge 0}\quad {\text resp.,}\quad
Y\times \A^1_K \stackrel{/\Z_d}{\longrightarrow}
\Spec\, (A_f)_{\le 0}\,.$$ \elem

The following two constructions will be important
in our applications below.

\bexa\label{EXA1} (Weighted blowup)
Let $S$ be an indeterminate of degree $-1$.
Consider the
graded subring
$$
\hA:=A[S]_{\ge 0}\cong \bigoplus_{\nu\ge \mu}
A_\nu S^\mu
$$
of the ring $A[S]$. By definition,
the {\it weighted blowup} of $V=\Spec A$
is the scheme $V':=\Proj \hA$.
Note that for every element $fS^\mu\in \hA$
with
$f\in A_{\nu}$ and $\nu>\mu$ we can write
$(fS^\mu)^n= f\cdot
(f^{n-1}S^{\mu n})$, where
$\deg f^{n-1}S^{\mu n}\ge 0$ for
$n\gg 0$. Therefore
$D_+(fS^\mu \hA)\subseteq D_+(f \hA)$, and so
$V'$ is covered by the affine  open subsets
$U_f:=D_+(f \hA)= \Spec \hA_{(f)}$, where
$$
\hA_{(f)}=\bigoplus_{\nu\ge 0}
\left( A_f\right)_{\nu}S^{\nu}
\cong (A_f)_{\ge 0}\cong F[T]^{\Z_{d}}
$$
with $T$ being an indeterminate of degree $1$.
Thus, $U_f\cong (Y\times \A^1_K)/\Z_d$,
see \ref{il} (b). In particular, if $K_Y$ is
a $\Q$-Cartier divisor then so is $K_{U_f}$,
see \ref{fingr} (c).
Notice also that the blowup morphism
$\sigma: V' \to V=\Spec A$ restricted to $U_f$
is induced by the inclusion
$A\hookrightarrow (A_f)_{\ge 0}$, and that
the exceptional divisor $E=\sigma^{-1}(V(A_+))$
is isomorphic
to $\Proj A$ under the natural morphism
$V'=\Proj \hA\to \Proj A$.
\eexa

\bexa\label{EXA2} (Weighted completion)\ \
Let $T$  be an indeterminate  with $\deg\,T=1$,
and consider
the projective scheme
$\bar V=\Proj  A[T]$. The scheme $V=\Spec A$
is naturally isomorphic  to the
affine open subset $D_+(T)\subseteq \bar V$ as
$$
A[T]_{(T)}=\bigoplus_{\nu\ge 0}
A_{\nu}T^{-\nu}\cong A\,.
$$ Thus $\bar V=V\cup D_{\infty}$ (where
$D_{\infty}:=\{T=0\}\cong \Proj A$
is the  `divisor at infinity')
is indeed proper over $\Spec A_0$; we call $\bar V$ {\it
the weighted completion} of $V$. The divisor
$D_{\infty}$ is covered by the affine open subsets
$U^f:=D_+(fA[T])\simeq \Spec
(A_f)_{\le 0}$ with
$f\in A_{d},\,\,d>0$ (indeed,
$$
A[T]_{(f)}=(A_f[T])_0=\bigoplus_{\nu\ge 0}
(A_f)_{-\nu}T^{\nu}
\cong (A_f)_{\le 0}\,{\rm ).}
$$ Furthermore, by \ref{il} (b)
we have
$(A_f)_{\le 0}\cong F[S]^{\Z_{d}}$,
where this time $S=T^{-1}$ is an indeterminate
of degree
$-1$, and again $U^f\simeq (Y\times \A^1_K)/\Z_d$.
\eexa

\brem If,
more generally,  $A$ is a graded (but not
necessarily
positively  graded) ring,
then  the above construction
provides a partial completion $\bar V$ of $V$.
\erem

   From \ref{a.1}-\ref{il} we obtain the
following proposition.

\bprop\label{quotsing}
With the assumptions
of \ref{a.1},
if $V^*=V\backslash V(A_+)$
is smooth then $\Proj A$
as well as the weighted blowup
$\Proj A[S]_{\ge 0}$, $\deg S=-1$,
have at most cyclic
quotient singularities.
Similarly, the subset
$\bar V\backslash V_+(A_+A[T])$ of the weighted
completion
$\bar V=\Proj A[T]$, $\deg T=1$,
has at most cyclic quotient
singularities.
\eprop

\proof
For a homogeneous element $f\in A_d$, $d>0$, the
morphism
$$\ul{\mu}:Y\times \G_m\to \Spec A_f
=V\backslash D_f\subseteq V^*$$
is \'etale (see \ref{il} (a)),
and so $Y=\Spec F$ is smooth.
Hence
$\Spec A_{(f)}\cong \Spec
\left(F^{\Z_d}\right)=Y/\Z_d$
has at most cyclic quotient singularities.
As
$\Proj A$ is covered by the affine open
subsets $D_+(f)=\Spec
A_{(f)}$, this scheme has
at most cyclic
quotient singularities as well.
The
proof of the remaining cases is
similarly
using the descriptions of
affine open coverings given in \ref{EXA1}
and \ref{EXA2}.
\qed\medskip

\subsection{Characterizing log
canonical forms in
terms of gradings}
We fix the following notation.

\bnota \label{notation}
Let now $K=\C$, and let
$A$ be a normal $\C$-algebra of
finite type with a
grading $A=\bigoplus_{\nu\ge 0} A_{\nu}$
(which
corresponds to a $\C^*$-action on $A$ via
$t.a=t^{\deg a}a$, where $t\in \C^*$
and $a\in A$ is homogeneous).
In the sequel such an affine variety
$V=\Spec A$  with
$\C^*$-action  will be referred to
as a {\it quasihomogeneous variety}
(indeed, there is a closed affine
embedding
$V\stackrel{i}{\hookrightarrow}\C^n$
given by a set of homogeneous generators
of $A$ and
equivariant
with respect to a diagonal $\C^*$-action on
$\C^n$).
Note that the {\it vertex set} $V(A_+)$
is just the fixed point set of the $\C^*$-action
on $V$.

Fix a $\C^*$-invariant divisor $D$ on $V$.
In the next theorem we
compute the $A$-modules
$$
L^{m}_{A,D}:=H^0(V,\cL^{m}_{V,D})
\quad\mbox{and}\quad
L^{2,m}_{A,D}:=H^0(V,\cL^{2,m}_{V,D})
$$
as submodules of
$$
H^0(V^*,\cL^{m}_{V,D})
\quad\text{respectively}\quad
H^0(V^*,\cL^{2,m}_{V,D}),
$$
where as before $V^*=V\backslash V(A_+)$.
Note
that all these modules have an induced
$\C^*$-action  and so they carry natural gradings.
\enota

\bthm\label{comp} If $V=\Spec
A$  with
$A=\bigoplus_{\nu\ge 0} A_{\nu}$
is a normal quasihomogeneous variety and
$D\subseteq V$ is a $\C^*$-invariant divisor
containing the divisorial part of $V(A_+)$ then
\be\label{coinc}
    \qquad L_{A,D}^{m}=H^0(V^*,\cL^{m}_{V,D})_{\ge 0}
\quad {\rm and}\quad
L_{A,D}^{2,m}=H^0(V^*,\cL^{2,m}_{V,D})_{>0}\,.\ee
\ethm

\bsit \label{sc} Observe first that if
$B:=(A_f)_{\ge 0}$ with $f\in A_+$ homogeneous,
then $(B_f)_{\ge 0}=B$. Bearing this in mind
we start with the following special case.
\esit

\blem\label{complem}
With the notation and assumptions as in
\ref{comp},
suppose furthermore
that there is an element $f\in A_d$, $d>0$,
such that
$A=(A_f)_{\ge 0}$. Then $V^*=\Spec A_f$, and (\ref{coinc})
holds, i.e.\
$$
L_{A,D}^{m}=(L_{A_f,D^*}^{m})_{\ge 0}
\qquad {\rm and}\qquad
L_{A,D}^{2,m}=(L_{A_f,D^*}^{2,m})_{>0}\,
$$
with $D^*:=D\cap V^*$.
\elem

\proof
As in \ref{il}, consider  the homogeneous
homomorphism
$$
\mu: A= (A_f)_{\ge 0}\longrightarrow F[T]
\quad\text{with}\quad
a/f^k \longmapsto
\bar a \cdot T^{\deg a-dk},
$$
where $T$ is an indeterminate with
$\deg\,T=1$ and ${\bar a}$ is
the residue class of $a$ in $F=A/(f-1)$.
By \ref{il} this
homomorphism induces an \'etale morphism
$A_f\to F[T,T^{-1}]$, and
$$
A_f\cong F[T,T^{-1}]^{\Z_d} \,,\quad
A=(A_f)_{\ge 0}\cong F[T]^{\Z_d}\, ,
\quad (A_f)_+\cong TF[T]^{\Z_d}\,
$$
with $\Z_d:=\Z/d\Z$.
Geometrically this means that, with
$Y=\Spec F$, the
morphism
$$
\ul{\mu}: Y\times \C \to V=\Spec A
$$
gives a  cyclic covering
(whence $V\cong (Y\times \C)/\Z_d$) non-ramified
off $\ul{\mu}^{-1}(V(A_+))=Y\times\{0\}$,
so that the restriction
$$\ul{\mu}|_{Y\times \C^*} : Y\times \C^*
\to V^*=V\backslash V(A_+)=\Spec A_f$$
is an \'etale covering.

To compute $L^{m}_{A,D}$ as a graded
submodule of  $L^{m}_{A_f,D^*}$,
notice first that by \ref{products} (a)
$$
L^{m}_{F[T],\tD}
\cong
L^{m}_{F,D\cap Y}\otimes L^{m}_{\C[T],\{0\}},
$$
where
$$\tD:=\ul{\mu}^{-1}(D)=(D\cap Y)
\times\C\cup Y\times\{0\}
\subseteq Y\times\C$$
(note that by our assumptions
$V(A_+)$ is a divisor contained in $D$).
The module $L^{m}_{\C[T],\{0\}}$ is equal to
$\C[T]\cdot (dT/T)^{\otimes m}$
and so it embeds into $L_{\C[T,T^{-1}]}^m$
as the submodule
of elements of degree $\ge 0$. It
follows that
$$L^{m}_{F[T],\tD}= (L^{m}_{F[T,T^{-1}],
(D\cap Y)\times\C^*})_{\ge 0}\,.$$
Taking invariants with respect to $\Z_d$
and using \ref{fingr} (a) we obtain that
$L^{m}_{A,D}=(L^{m}_{A_f,D^*})_{\ge 0}$.
The proof of the equality
$L^{2,m}_{A,D}=(L^{2,m}_{A_f,D^*})_{>0}$
is similarly using the fact
that $L^{2,m}_{\C[T],\{0\}}$ is generated by
$T(dT/T)^{\otimes m}$.
\qed\medskip

\medskip\no{\it Proof of Theorem \ref{comp}}.\,
Let $S$ be an indeterminate with
$\deg\,S=-1$.  Consider the weighted blowup
(see \ref{EXA1})
$$
\sigma:V':=\Proj (A[S]_{\ge 0})\to V\,
$$
with exceptional divisor $E\subseteq V'$,
and denote $D':=\sigma^{-1}(D)_\red
\cup E\subseteq V'$.
According to \ref{dirim} (c)
$H^0(V',\cL^{m}_{V',D'})\cong L^{m}_{A,D}$,
and so we
need to show that
$$
H^0(V',\cL^{m}_{V',D'})\cong
H^0(V^*,\cL^{m}_{V,D})_{\ge 0}\,.
$$
By \ref{EXA1} the affine open subsets
$$
U_f=D_+(f\hA)=\Spec(A_f)_{\ge 0} \subseteq V'
$$
with $f\in A_+$ homogeneous
form a covering of $V'$. If
$\omega\in H^0(V^*,\cL^{m}_{V,D})$ is a
homogeneous form, then
$\sigma^*(\omega)$ belongs to
$H^0(V',\cL^{m}_{V',D'})$ if
and only if, for all $f$
as above, the form $\omega|_{D(f)}\in
L^m_{A_f, D\cap D(f)}$
(where $D(f)=\Spec A_f\subset V^*$,
see \ref{sl})
extends to a form in
$L_{(A_f)_{\ge 0},D\cap U_f}^{m}$. As by
\ref{sc} and \ref{complem}
$$L_{(A_f)_{\ge 0},D\cap U_f}^{m}=
\left(L_{A_f,D\cap D(f)}^{m}\right)_{\ge 0}\,,$$
the result for $L^m$-forms follows.
The proof in the
case of
$L^{2,m}$-forms is analogously and left to
the reader.
\qed\medskip

\begin{sit}\label{sit}
In the next proposition we give a
dual version of \ref{comp} which allows
to control the $\cL^m$ and
$\cL^{2,m}$-forms at infinity
(this will be useful later on,
see \ref{pluqua0}).
With the
assumptions as in \ref{comp},
consider the weighted
completion $\bar V=\Proj A[T]$ with
$\deg T=1$ (see \ref{EXA2}).
The subset
$V_{\infty}:=\bar V\backslash V_+(A_+A[T])$
contains the  divisor at infinity
$D_\infty\cong \Proj A$
of $\bar V$. Let $\bar D$
denote the union $D\cup D_\infty$.
\end{sit}

\bprop\label{infin}
With the above notation we have
$$
\begin{array}{l}
H^0(V_{\infty},\,\cL_{V_{\infty},\bar D}^{m})=
H^0(V^*,\cL^{m}_{V,D})_{\le 0} \\[3pt]
H^0(V_{\infty},\,\cL_{V_{\infty},\bar D}^{2,m})=
H^0(V^*,\cL^{2,m}_{V,D})_{< 0}\,.
\end{array}
$$
\eprop

\proof
The variety $V_{\infty}$
is covered by the affine open
subsets $U^f=\Spec\, (A_f)_{\le 0}$,
where $f\in A_+$ is
homogeneous (see \ref{EXA2}). Applying \ref{complem} to
$(A_f)_{\le 0 }$ with the grading
reversed we obtain
$$
L_{(A_f)_{\le 0 },\bar D}^{m}
\cong (L_{A_f,\bar D}^{m})_{\le 0}
\quad\text{and}\quad
L_{(A_f)_{\le 0 },\bar D}^{2,m}
\cong (L_{A_f,\bar D}^{2,m})_{<0}\,.
$$
Now we can proceed as in the proof
of \ref{comp}; we leave the details
to the  reader.\qed\medskip

\subsection{Log terminal and log canonical
singularities of quasihomogeneous varieties}
With the notations as in \ref{notation}, let
$\omega_A:=H^0(V,\cO_V(K_V))$ be the dualizing module of
$A$ and let
$$
\omega_A^{[m]}(kD):=H^0(V,\cO_V(mK_V+kD)).
$$
We have seen in \ref{divisorial part} (b) that $L^m_{A,D}$
and $L^{2,m}_{A,D}$ are in a natural way submodules of
$\omega_A^{[m]}(mD)$ resp.\
$\omega_A^{[m]}((m-1)D)$. In the next theorem we
identify the modules
$$\gL_{A,D}^{m}:= H^0(V,\gL_{V,D}^{m})
\quad {\text and}\quad
\gD_{A,D}^{m}:=H^0(V,\gD_{V,D}^{m})$$
with certain graded pieces of $\omega_A^{[m]}(mD)$ resp.\
$\omega_A^{[m]}((m-1)D)$.

\bthm\label{comprem}
With the notation and assumptions
as in \ref{comp}
the following hold.

\begin{enumerate}
\item[(a)] If $(V^*, D^*)$ is lc then
$$
L_{A,D}^{m}\cong \omega_A^{[m]}(mD)_{\ge 0}
\quad {\text and}\quad
\gL_{A,D}^{m}\cong \omega_A^{[m]}(mD)_{< 0}\,.
$$
\item[(b)] Similarly, if  $(V^*, D^*)$ is
$L^2$-lt then
$$
L_{A,D}^{2,m}\cong \omega_A^{[m]}((m-1)D)_{> 0}
\quad {\text and}\quad
\gD_{A,D}^{m}\cong \omega_A^{[m]}((m-1)D)_{\le 0}\,.
$$
\item[(c)] If $(V^*,D^*)$ is
lc, respectively $L^2$-lt
then so is the weighted blowup $(V',D')$ (see
\ref{EXA1}).
\end{enumerate}
\ethm

\proof
(a)
Observe that
in virtue of \ref{Kempf} (b) $\gL^m_{V^*,D^*}=0$ or,
equivalently,
$\cL^m_{V^*,D^*}=\cO_{V^*}(m(K_{V^*}+D^*))$.
Hence
\ref{comp} implies that
$$
L_{A,D}^{m}\cong H^0(V^*,\cO_{V}(m(K_{V}+D)))_{\ge 0}.
$$
The module on the right contains
$H^0(V,\cO_{V}(m(K_{V}+D)))_{\ge 0}=
\omega_A^{[m]}(mD)_{\ge
0}$, and
by \ref{divisorial part}
(b) $\omega_A^{[m]}(mD)$
contains $L_{A,D}^{m}$ whence
$(\omega_A^{[m]}(mD))_{\ge 0}$ contains
$L_{A,D}^{m}=(L_{A,D}^{m})_{\ge 0}$. Thus
$L_{A,D}^{m}$  and $\omega_A^{[m]}(mD)_{\ge
0}$ are equal, as required. The proof of
(b) is similarly and left to the reader.

In order to show (c)  consider
first the special case treated in
\ref{complem} so that $A=(A_f)_{\ge 0}$. Then
$U_f=D_+(f\hA)=\Spec\, (A_f)_{\ge 0}   = V$ and
$\sigma |_{U_f}={\rm id}$, whence
$V'=V=(Y\times\C)/\Z_d$ (see \ref{il} (b)). If
$(V^*, D^*)$ is lc then (with the notation
as in the proof of \ref{complem}) the
\'etale covering
$$
(Y\times\C^*,\ul{\mu}^{-1}(D\cap V^*))=
(Y\times\C^*,(D\cap Y)\times\C^*)
$$
of $(V^*, D^*)$ is also lc.
By \ref{products} $(Y,D\cap Y)$ is lc and so,
applying \ref{products} again,
$(Y\times\C,\tD)$
(with $\tD=(D\cap Y)\times\C
\cup Y\times\{0\}=
\ul\mu^{-1}(D)\subseteq Y\times \C$)
is lc as well.
As $\ul\mu |_{(Y\times\C^*)}$ is unramified,
the divisor $\ul\mu^*(K_V+D)$ on $Y\times\C$
is
equal to $K_{Y\times\C}+\tD+
\lambda (Y\times\{0\})$
for some $\lambda \in\Z$;
in particular, it is
$\Q$-Cartier. Taking quotients and using
\ref{fingr} and \ref{Kempf} (b) we deduce that
$(V, D)=(V', D')$ is also lc.

The general case follows easily
from this with the same
reasoning as in the proof of \ref{comp}.
The proof in the $L^2$-lt case is similarly
and left to the reader.
\qed\medskip

\bcor\label{cor logcan}
With the assumptions as in \ref{comp},
the following hold.

\begin{enumerate}
\item[(a)] If $(V^*,D^*)$ is lc and $K_V+D$ is
$\Q$-Cartier then
$(V,D)$ is lc if and only if
$$\gL^m_{A,D}=\omega_A^{[m]}(mD)_{<0}
=0\quad\forall m\ge 1\,.$$
\item[(b)] If $(V^*,D^*)$ is $L^2$-lt then
$(V,D)$ is $L^2$-lt if and only if
$$\gD^m_{A,D}=\omega_A^{[m]}((m-1)D)_{\le 0}=0
\quad\forall m\ge 1\,.$$
\item[(c)] {\em (\cite{Fl})}\ If $V^*$ has
rational singularities and $V$ is
Cohen-Macaulay then
$V$ has rational singularities if and only if
$(\omega_A)_{\le 0}=0$.
\end{enumerate}
\ecor

\proof As the coherent sheaves
$\gL^m_{V,D}$ resp., $\gD^m_{V,D}$
on the affine variety $V$
are globally generated, under
the assumptions of (a) resp.,
(b) they vanish. Now
(a) and (b) follow immediately
from \ref{comprem} (a), (b) and
\ref{Kempf} (b)-\ref{L2}. To prove (c)
observe that
$\omega_A=H^0(V, \cO_V(K_V))$.
Thus
(c) is a consequence of \ref{comprem} (b)
and \ref{Kempf} (a).
\qed\medskip

\bsit\label{Gor}
Assume now that $\omega_A$ is a free $A$-module so that
$$
\omega_A\cong A[N_A]
$$
for some $N_A\in\Z$ (where as usual
$A[N]$ is the module $A$ equipped
with the new grading $A[N]_i:=A_{i+N}$). Note that by the
homogeneous version of Nakayama's lemma this assumption is
satisfied if, for instance,
$A$ is a Gorenstein ring and $A_0$
is a local ring. If moreover
$D$ is given by a homogeneous equation of degree $d$
then
$$
\begin{array}{l}
H^0(V,\cO_{V}(m(K_{V}+D)))\cong
A[m(N_A+d)]\qquad {\rm and}\\[2pt]
H^0(V,\cO_{V}(mK_{V}+(m-1)D))\cong A[mN_A+(m-1)d]\, .
\end{array}
$$
  From \ref{cor logcan} we obtain
the following characterizations.
\esit

\bcor\label{corlogcan}
Let $(A,D)$
be as in \ref{comp}.
If
$\omega_A\cong A[N_A]$ and
$D$ is given by a homogeneous
equation of degree $d$
then the following hold.
\begin{itemize}
\item[(a)] If $(V^*,D   ^*)$ is lc then
$(V,D)$ is lc if and only if $N_A+d\le 0$.
\item[(b)] If $(V^*,D   ^*)$
is $L^2$-lt then
$(V,D)$ is $L^2$-lt if and only if $d>0$
and $N_A+d\le 0$ or $d=0$ and $N_A<0$.
\item[(c)] \cite{Fl} If $V^*$
has rational singularities and $A$ is
Cohen-Macaulay then $V$
has rational singularities
if and only if $N_A< 0$.
\end{itemize}\ecor

\bexa\label{ci logcan}
In particular,
let $V=\Spec A$ be a normal
complete intersection
of  dimension
$n$ given in $\C^{n+s}$ by polynomials
$p_1,\dots, p_{s}\in \C^{n+s}$
which are  quasihomogeneous of
degrees $d_1,\dots, d_s$
with respect to
weights $w_j\ge 0,\,\,j=1,\dots,n+s$, so that
$$
p_i(\lambda^{w_{1}}x_1,\dots,
\lambda^{w_{n+s}}x_{n+s})=
\lambda^{d_i}p(x_1,\dots,x_{n+s}),\qquad
i=1,\dots,s\, .
$$
It is well known (see e.g.\ \cite[p.42]{Fl})
that
$A$ is Gorenstein (whence Cohen-Macaulay)
and
\be\label{loca}
\omega_A=A[N_A]\quad\text{with}\quad
N_A (=N_V):=\sum_{i=1}^{s} d_i -
\sum_{j=1}^{n+s}w_{j}\, .
\ee
\eexa

\bexa\label{EVA} As a concrete example,
consider the polynomial
$x^ay^d+u^b+v^c=0$ with $a,b,c,d\ge 2$.
It is weighted homogeneous of
degree $abc$ with weights, say,
$\deg x=bc$, $\deg y=0$, $\deg v=ab$
and $\deg u=ac$.
The associated
quasihomogeneous variety
$V=\{p=0\}\subseteq\C^4$
has singularities only along the
$x$- and $y$-axes.
By \ref{corlogcan}
the singularities along the $x$-axis off the
origin
are log canonical iff $1/b+1/c+1/d\ge 1$;
they are rational and, moreover, log
terminal iff $1/b+1/c+1/d> 1$ (indeed,
locally near a point $(x_0,0,0,0)\in V^*$
with $x_0\neq 0$
the mapping $$V\ni (x,y,u,v)\longmapsto
(x, x^{a/d}y,u,v)\in {\tilde V}:=\{y^d+u^b+v^c=0\}
\subseteq\C^4$$ is well defined and biholomorphic).
Letting $D$ be the divisor $\{x=0\}$ it follows that
$(V,D)$ is lc if and only if
$1/b+1/c\ge 1$.
Moreover the singularities of $V$
are rational (and log terminal)
if and only if $1/a+1/b+1/c> 1$
and $1/b+1/c+1/d> 1$.
\eexa

As another example we study varieties
that are given
by the maximal Pfaffians of a
skew symmetric matrix (cf.\ \cite{BEi}).

\bexa\la{pfaff}
Let $R=\bigoplus_{i\ge 0} R_i$
be a finitely generated graded
$\C$-algebra which is Gorenstein
with $\omega_R\cong R[N_R]$.
Consider a skew symmetric
$(2n+1)\times (2n+1)$-matrix $(a_{ij})$
of homogeneous elements of
$R$ and assume that $\deg
a_{ij}=d_i+d_j-N$, where
$N$ and $d_i$, $1\le i\le 2n+1$,
are positive integers
satisfying $\sum_id_i=nN$.
The maximal  Pfaffians generate a
homogeneous ideal, say, $I$ of $R$.
Assume that the
quotient $A:=R/I$ has dimension $\dim R-3$.
Then by
\cite{BEi} the minimal resolution of $A$
has the form
$$
0\to R[-N]\to \bigoplus_{j=1}^{2n+1}
R[-N+d_j]\stackrel{(a_{ij})}{\lto}
\bigoplus_{i=1}^{2n+1} R[-d_i] \to R \to
A\to 0\,.
$$
Taking $\Ext_R(-, \omega_R)$ of this sequence
it follows that
$\omega_A\cong A[N_A]$ with $N_A:=N+N_R$.
Hence
\ref{corlogcan} applies in this situation.
\eexa

\subsection{Plurigenera of quasihomogeneous
singularities}

The results above provide the following
explicit formulas for the plurigenera of
isolated singularities of
quasihomogeneous varieties, where as usual
$\omega_A^{[m]}$
denotes the reflexive hull of the module
$\omega_A^{\otimes m}$
of K\"ahlerian $m$-differentials on $A$.

\bprop\label{del} Let
$A=\bigoplus_{\nu\ge 0} A_{\nu}$ be a normal
$\C$-algebra of
finite type and assume that the corresponding
affine variety $V=\Spec A$
has $\dim V \ge 2$ and
at most isolated singularities. Then the
following hold.

\begin{enumerate}
\item[(a)] If
$A_0\ne \C$ (in particular,
if $V$ has at least two singular points)
then  $\delta_m (V,p)=\lambda_m (V,p)=0$
for all $m\ge 1$ and $p\in\Sing V$.
\item[(b)] If $A_0=\C$ and
$V$ has a  unique singular point
$p$ then for every $m\ge 1$ we have
$$
\delta_m (V,p)=\dim_{\C}\,(\omega_A^{[m]})_{\le
0}\qquad {\rm and}\qquad
\lambda_m (V,p)=\dim_{\C}\,(\omega_A^{[m]})_{<0}\,.
$$
\end{enumerate}
\eprop

\proof (a) Let $D$ denote the divisorial part of
$V(A_+)$. Using
\ref{injective} resp.\  \ref{comprem} (b) we have an
inclusion
$$
\gD_A^m\subseteq\gD_{A,D}^m\cong
(\omega_A^{[m]})((m-1)D)_{\le 0}.
$$
By our
assumption, $\dim_\C A_0=\infty$.
Since
$(\omega_A^{[m]})((m-1)D)$ is
a torsion-free $A_0$-module, its $A_0$-submodule
$\gD_A^m$ is also
torsion-free over $A_0$, and so either $\gD_A^m=0$ or
$\dim_{\C} \gD_A^m=\infty$. On the other hand,
$$
\dim_{\C}\,\gD^m_A=\sum_{p\in\Sing V}
\delta_m (V,p)
<\infty\,.
$$
It follows that $\gD_A^m=0$,
and so $\delta_m (V,p)$
vanishes for all $m\ge 1$ and
$p\in\Sing V$. Consequently
$\lambda_m (V,p)=0$, and (a)
follows.

(b) Since by our assumption $V_+=\Spec A_0=\{p\}$
and so ${\rm div}\, V_+=\emptyset$,
by \ref{comprem} (a),
(b) (with $D=0$) we have
$$
\gL_A^m\cong (\omega_A^{[m]})_{< 0}
\qquad {\rm and}\qquad
\gD_A^m\cong (\omega_A^{[m]})_{\le 0}\,.
$$
Therefore (since $p$ is the unique singular
point of $V$) we obtain
$$
\delta_m (V,p)=\sum_{x\in\Sing V}
\delta_m (V,x)=
\dim_{\C}\,\gD^m_A =
\dim_{\C}\,(\omega_A^{[m]})_{\le 0}\,,
$$
and similarly
$\lambda_m (V,p)=\dim_{\C}\,(\omega_A^{[m]})_{<0}$,
proving (b).
\qed\medskip

\bexa\label{ci}
Let $V=\Spec A$ be a complete intersection
as in \ref{ci logcan} and
assume moreover that
$V$ has an isolated singularity at the origin
$0\in \C^{n+s}$. As $\omega_A=A[N_A]$  we have
$\omega_A^{[m]}=A[mN_A]$.  Now \ref{del} (a)
implies the result of \cite{Mo} which says that
$$ \delta_m=\dim\,
\sum_{i\le 0} A_{i+mN_A}\qquad
{\rm and}\qquad
\lambda_m=\dim\,\sum_{i< 0} A_{i+mN_A}\,.$$
\eexa

\subsection{Log-Kodaira dimension
of quasihomogeneous varieties}
Recall the following notions.

\bdefi\label{lkd} {\rm \cite{Ii}
Let  $V$ be a
smooth quasi-projective variety, and
let $\bar V$ be a
smooth compactification  of $V$ by an
SNC divisor
$\bar D=\bar V\backslash V$. The {\it
logarithmic
plurigenera}
$\bar p_m(V)$ are defined by
$$\bar p_m(V)=\dim\,H^0(\bar
V,\,\cO_{\bar V}(m(K_{\bar
V}+\bar D)))\,,\,\,\,m\ge 1.$$
The {\it logarithmic Kodaira dimension}
$\bar k\,(V)$ of $V$ is
$$
\bar k\,(V)=
\left\{ \begin{array}{cc}
-\infty & {\rm if}\quad \bar p_m(V)=0
\quad\forall
m\in\N\\
\min\,\{k\in \N :  \limsup\limits_{m\to\infty}\bar
p_m(V)/m^k<\infty\} & \quad{\rm otherwise.}
\end{array} \right .
$$}
\edefi

Taking the intersection of the modules in
\ref{comp}
(or \ref{comprem}) and
\ref{infin}  we obtain the following
proposition.

\bprop\label{pluqua0}
Let $A=\bigoplus_{i\ge 0}A_i$ be
a normal graded $\C$-algebra
of finite type with
$A_0=\C$
(so that $V(A_+)=\{p\}$)
and let $V=\Spec A$
be the corresponding
quasihomogeneous variety.
If $D$
is a homogeneous
reduced divisor on V then
$$
\bar p_m(V\backslash (D\cup \Sing V))
=\dim\,(L^{m}_{A,D})_0\,.
$$
In particular, $\bar k(V\backslash(D\cup
\Sing V))=-\infty$
if and only if
$(L^{m}_{A,D})_0=0$ for all $m\ge 1$.
\eprop

Summarizing the preceding results gives the following
theorem.

\bthm\label{pluqua}
With the assumptions as in \ref{pluqua0}
above, suppose in addition
that the divisor $D$
is given by a
homogeneous equation of degree $d$
and that the pair $(V^*,D^*)$ is lc
(where as before,
$V^*:=V\backslash \{p\}$  and $D^*=D\backslash \{p\}$
with
$\{p\}=V(A_+)$ being the vertex of $V$).
Then the following hold.

\begin{enumerate}
\item[(a)]
$\bar p_m(V\backslash (D\cup \Sing V))=
\dim\,(\omega_A^{[m]})_{md}\,.$

\item[(b)] If moreover $V$ has an isolated
singularity at $p$, then
$\bar k(V^*)=-\infty$ if and only if
$(\omega_A^{[m]})_0=0$ for all $m\ge 1$ or,
equivalently,
if and only if $\delta_m(V,\,p)=0$ for all $m\ge 1$.

\item[(c)] Suppose that $K_V$ is a $\Q$-Cartier
divisor,
so that for a certain $m_0\ge 1$ we have
$\omega_A^{[m_0]}\cong A[N]$  with some
$N\in\Z$. Then
$$
\bar k(V\backslash(D\cup \Sing V))=
\left\{ \begin{array}{ccc}
-\infty & {\rm if}\quad N+m_0 d<0\\
0 & {\rm if}\quad N+m_0 d=0\\
\dim V -1 & {\rm if}\quad N+m_0 d>0\,.
\end{array} \right .
$$
\item[(d)] In particular,
with the assumptions as in (c)
$\bar k(V\backslash(D\cup \Sing V))\le 0$
if and
only if $(V,D)$ has a log canonical singularity at
$p$. If moreover $(V,p)$ is an
isolated singularity
then $(V,p)$ is lt
if and only if
$\bar k(V^*)=-\infty$.
\end{enumerate}
\ethm

\proof
By \ref{comprem} (a) we have
$L^{m}_{A,D}\cong
\omega_A^{[m]}[md]_{\ge 0}$
and so (a)
follows from \ref{pluqua0}.
Under the assumptions as in (b),
by \ref{pluqua0} we have
$p_m(V^*)=(\omega_A^{[m]})_0$.
Thus the first equivalence
in (b) is an immediate consequence of (a).
To show the second equivalence,
choose a non-zero homogeneous element $g\in A$ of
some degree, say $k>0$. If
$\eta$ is a non-zero form of degree $s<0$ in
$\omega^{[m]}_A$, then $g^s\eta^k$ is a
non-vanishing form of degree 0 in
$\omega^{[mk]}_A$. In other words,
$(\omega^{[m]}_A)_0 = 0$ for all $m\ge 1$ if and
only if $(\omega^{[m]}_A)_{\le 0}=0$ for all $m\ge
1$. In view of \ref{del} (a) this proves (b).

In order to deduce (c) notice that (in virtue of
(a))
$$
\bar p_m(V\backslash (D\cup \Sing V))=
(\omega_A^{[km_0]})_{km_0 d}=A_{k(N+m_0 d)}=0\quad
\mbox{for all } k\ge 1
$$
if and only if $N+m_0 d<0$.
In the case $N+m_0 d=0$ we have
$(\omega_A^{[km_0]})_{km_0
d}=A_0= \C$ for all $k\ge 1$.
If $N+m_0 d>0$ then
we may choose a number
$r$ such that
$\dim\,(\omega_A^{[rkm_0]})_{rkm_0 d}
=\dim\,A_{rk(N+m_0d)}$
grows asymptotically like
$k^{\dim A-1}=k^{\dim \Proj A}$. (d) follows from
\ref{kmp} and
\ref{corlogcan} (a) above.
\qed\medskip

\brem\label{vertex pos rem}
Let $V=\Spec A$ be a normal quasihomogeneous
variety with a vertex set
$V_0=\Spec A_0$ of positive dimension
(that is, $A_0\neq \C$). If $V$
has
isolated singularities then
$\bar k(V\backslash \Sing V)=-\infty$
(cf. \ref{del} (b)).
Indeed, the general fibre of the canonical
projection
$q:V\to V_0$ is  a smooth quasihomogeneous
variety
with a vertex set of
dimension $0$, whence it is isomorphic
to $\C^k$ for some $k>0$
(see e.g., \cite[8.5]{Za}).
\erem

\bcor\label{neco}
Let $A:=\C[X_1,\dots,X_{n+s}]/(p_1,\dots,p_s)$
be a normal
quasi\-homo\-geneous complete intersection with
    weights
$w_1,\dots,w_{n+s}>0$ (see \ref{ci logcan}).
With $V=\Spec
A$ being the associated affine variety and $D$
being a homogeneous degree
$d$ divisor on $V$, we have
$$
\bar k(V\backslash(D\cup \Sing V))=
\left\{ \begin{array}{ccc}
-\infty & {\rm if}\quad N_A+d<0\\
0 & {\rm if}\quad N_A+d=0\\
\dim V -1 & {\rm if}\quad N_A+d>0\,.
\end{array} \right .
$$
\ecor

\proof
As the weights are all positive we have
$A_0=\C$, moreover
$\omega_A=A[N_A]$ (see (\ref{loca})
in \ref{ci logcan}),
and so \ref{pluqua} (c) (with $m_0=1$)
applies.  \qed
\end{document}